\documentstyle[amssymb,amsfonts,12pt]{amsart}
\newtheorem{theorem}{Theorem}[section]

\newtheorem{corollary}[theorem]{Corollary}
\newtheorem{proposition}[theorem]{Proposition}

\newtheorem{conjecture}[theorem]{Conjecture}

\theoremstyle{remark}
\newtheorem{remark}[theorem]{Remark}

\def\QSet{\mbox{\rm\kern.24em
\vrule width.03em height1.48ex depth-.051ex \kern-.26em Q}}

\def\R{{\mathbb R}}

\def\N{{\mathbb N}}
\def\C{{\mathbb C}}

\def\Z{{\mathbb Z}}

\def\be#1{\begin{equation}\label{#1}}

\def\bas{\begin{align*}}
\def\eas{\end{align*}}
\def\bi{\begin{itemize}}
\def\ei{\end{itemize}}
\newenvironment{proof}{\noindent {\bf Proof} }{\endprf\par}
\def \endprf{\hfill  {\vrule height6pt width6pt depth0pt}\medskip}
\def\emph#1{{\it #1}}

\begin{document}

\title[Proof of the HRT conjecture for (2,2) configurations]{Proof of the HRT conjecture for (2,2) configurations}

\author{Ciprian Demeter}
\address{Department of Mathematics, Indiana University, 831 East 3rd St., Bloomington IN 47405}
\email{demeterc@@indiana.edu}

\author{Alexandru Zaharescu}
\address{Department of Mathematics, University of Illinois at Urbana-Champaign, 1409 W. Green Street, Urbana, Illinois 61801-2975}
\email{zaharesc@@math.uiuc.edu}

\keywords{}
\thanks{The  first author is supported by a Sloan Research Fellowship and by NSF Grants DMS-0742740 and 0901208}
\thanks{The  second author is supported by  NSF Grant DMS-0901621}
\thanks{ AMS subject classification: Primary 26A99; Secondary 11K70, 65Q20}
\begin{abstract}

We prove that for any 4 points in a (2-2) configuration, there is no linear dependence between the associated time-frequency translates of any  $L^2(\R)$ function.
\end{abstract}
\maketitle

\section{Introduction}
The following conjecture, known as the HRT conjecture appears in \cite{HRT}. See also \cite{H} for an ample discussion on the subject.

\begin{conjecture}
\label{cc1}
Let $(t_j,\xi_j)_{j=1}^n$ be $n\ge 2$ distinct points in the plane. Then there is no nontrivial $L^2$ function  $f:\R\to \C$ satisfying a  nontrivial linear dependence
$$\sum_{j=1}^nd_if(x+t_j)e^{2\pi i\xi_jx}=0,$$
for a.e. $x\in \R$.
\end{conjecture}

The conjecture  follows trivially when the points $(t_j,\xi_j)_{j=1}^n$ are collinear. The  conjecture was proved  when  $(t_i,\xi_i)_{i=1}^n$ sit on a lattice, \cite{Lin}, using von Neumann algebras techniques. See also \cite{BS}, \cite{DG}, for more elementary alternative arguments. In particular, this is the case with any 3 points. But the question whether the conjecture holds for {\em arbitrary} 4 points is open. Progress on that has been made by the first author in \cite{De2} using a number theoretical approach, and we briefly discuss it below.

We will call an $(2,2)$ configuration, any collection of 4 distinct points in the plane, such that there exist 2 distinct parallel lines each of which containing 2 of the points. One of the results in \cite{De2} is

\begin{theorem}
\label{TT2'}
Conjecture \ref{cc1} holds  for special $(2,2)$ configurations $(0,0),(1,0), (0,\alpha), (1,\beta)$

(a) if
$$\liminf_{n\to\infty}n\log n\min\{\|n\frac{\beta}{\alpha}\|,\|n\frac{\alpha}{\beta}\|\}<\infty$$

(b)  if at least one of $\alpha,\beta$ is rational

In either case, no nontrivial solution $f$ can exist satisfying  minimal decay
$$\lim_{|n|\to\infty\atop{n\in\Z}}|f(x+n)|=0,\;\;a.e.\;x$$
\end{theorem}

In this paper we prove the strongest possible statement about (2,2) configurations, namely

\begin{theorem}
\label{TT2}
Conjecture \ref{cc1} holds  for all $(2,2)$ configurations.
Moreover, when the points sit in a special (2,2) configuration $(0,0),(1,0), (0,\alpha), (1,\beta)$, no nontrivial solution $f$ can exist satisfying  minimal decay
$$\lim_{|n|\to\infty\atop{n\in\Z}}|f(x+n)|=0,\;\;a.e.\;x$$
\end{theorem}

The general approach for proving this theorem is the one developed in \cite{De2}. We first reduce to the case of special configurations, by applying metaplectic transformations.   Then we  turn the hypothetical linear dependence into a recurrence. The contribution from $\beta$ is estimated by using the {\em conjugates trick}. The novelty of our approach here is in the way we treat the contribution coming from the terms containing $\alpha$. In particular, we exploit the Diophantine behavior of $\alpha$ at more than one scale.

\section{Proof of the main theorem}
Define $[x]$, $\{x\}$, $\|x\|$  to be the integer part, the  fractional part and the distance to the nearest integer of $x$.  For two quantities $A$, $B$ that vary, we will denote by $A\lesssim B$ or $A=O(B)$ the fact that $A\le C B$ for some universal constant $C$, independent of $A$ and $B$.  In general, $A\lesssim_p B$ means that the implicit constant is allowed to depend on the parameter $p$. The notation $A\sim_p B$ means that $A\lesssim_p B$ and $B\lesssim _p A$.
 If no parameter is specified, the implicit constants are implicitly understood to depend on the (harmless) fundamental parameters introduced in the proof of Theorem \ref{TT2}. For a set $A\subset \R$, we will denote by $|A|$ its Lebesgue measure, and if the set is finite, $|A|$ will represent its cardinality. Finally, we define
$e(x):=e^{2\pi i x}$.

Let $0<\alpha<1$ be irrational.
Let $\frac{p_k}{N_k}$ be the $k^{th}$ convergent of $\alpha$, so that
\begin{equation}
\label{ecc6}
|\alpha-\frac{p_k}{N_k}|\le \frac1{N_kN_{k+1}},
\end{equation}
and
\begin{equation}
\label{ecc4}
p_{k}N_{k-1}-p_{k-1}N_{k}=(-1)^{k-1}
\end{equation}
Since
$$N_{k}\le N_{k+1},$$there exists an infinite set $E\subset \N$ and a constant $D=D(\alpha)$ such that for each $k\in E$ we have
\begin{equation}
\label{ecc1}
\frac{N_k}{N_{k+1}}\le D  \min_{j\le k}\frac{N_j}{N_{j+1}}.
\end{equation}
Define $\frac1{M_k}:=N_k^2|\alpha-\frac{p_k}{N_k}|$. Of course, $M_k\ge 1$ for each $k$.

The following proposition is the main new ingredient in this paper.

\begin{proposition}
\label{P1}
Let $k\in E$ be odd, and $0<\delta<\frac{1}{100}$. Define $N:=N_k$, $p:=p_k$, $M:=M_k$. Then, for each $x\in [0,1]$ such that \begin{equation}
\label{ecc2}
\min\{\frac{\|x\|}{N},\|x-n\alpha\|,\|x-\frac{n}{N}\|:1\le n\le N\}\ge \frac{\delta}{N}\end{equation} we have
\begin{equation}
\label{ecc36}
\prod_{n=1}^{N}|{e(x)-e(\alpha n)}|\sim_\delta 1.
\end{equation}
\end{proposition}
\begin{remark}
The key thing in \eqref{ecc36} is that the similarity constant does not depend on $N$.
\end{remark}
\begin{proof}
Fix $x$ satisfying \eqref{ecc2}. We will compare $\prod_{n=1}^{N}|{e(x)-e(\alpha n)}|$ to $$\prod_{n=1}^{N}|{e(x)-e(\frac{np}{N})}|=\prod_{n=1}^{N}|{e(x)-e(\frac{n}{N})}|=|e(Nx)-1|\sim_\delta 1,$$
and prove that their ratio is $\sim_\delta 1$.
This is reasonable to expect, since, due to \eqref{ecc6}, we have for each $1\le n\le N$
\begin{equation}
\label{ecc7}
|n\alpha-\frac{np}{N}|\le \frac{1}{N}
\end{equation}
First, let $1\le n_1,n_2,\ldots, n_{200}\le N$ be such that
$$\|x-\frac{n_ip}{N}\|\le \frac{100}N$$
Due to \eqref{ecc2} and \eqref{ecc7}, we get that
\begin{equation}
\label{ecc9}
\delta^3\lesssim \prod_{i=1}^{200}\frac{|{e(x)-e(\alpha n_i)}|}{|{e(x)-e(\frac{n_ip}{N})}|}\lesssim \delta^{-1}.
\end{equation}
Next, we analyze
$$\prod_{n=1\atop{n\not=n_i}}^{N}\frac{|{e(x)-e(\alpha n)}|}{|{e(x)-e(\frac{np}{N})}|}.$$
Note that
$$\frac{|{e(x)-e(\alpha n)}|}{|{e(x)-e(\frac{np}{N})}|}=\left|1+\frac{1-e(\alpha n-\frac{np}{N})}{e(x-\frac{np}{N})-1}\right|,$$
and that
$$\left|\frac{1-e(\alpha n-\frac{np}{N})}{e(x-\frac{np}{N})-1}\right|\le \frac{10}{N\|x-\frac{np}{N}\|}<\frac12.$$
Thus,
$$\sum_{n=1\atop{n\not=n_i\atop{\|x-\frac{np}{N}\|\ge \delta}}}^N\frac1{N\|x-\frac{np}{N}\|}\lesssim \sum_{N\delta\le i\le N}\frac{1}{i}\lesssim \log (\delta^{-1}).$$
Using this and the fact that
$$1+x\le e^x,\;\;0<x<1$$
$$e^{-10x}\le 1-x,\;\;0<x<1/2,$$
we get
\begin{equation}
\label{ecc3}
\prod_{n=1\atop{n\not=n_i\atop{\|x-\frac{np}{N}\|\ge \delta}}}^N\left|1+\frac{1-e(\alpha n-\frac{np}{N})}{e(x-\frac{np}{N})-1}\right|\sim_\delta 1
\end{equation}

Denote by
$$A:=\{1\le n\le N:n\not=n_i,\;\;\|x-\frac{np}{N}\|<\delta\}$$
Using the fact that for $z\in\R$ with $|z|<\frac1{10}$
$$1/2\le\frac{|e(z)-1|}{2\pi|z|}<2,$$we get for each $n\in A$
$$\left|\frac{1-e(\alpha n-\frac{np}{N})}{e(x-\frac{np}{N})-1}\right|< \frac{\frac{10n}{N^2M}}{\frac{100}{N}}<\frac{1}{10}.$$
It is easy to check that for each $z\in \C$ with $|z|<\frac1{10}$ we have
$$e^{-O(|z|^2)}\le\left|\frac{1+z}{e^z}\right|\le e^{O(|z|^2)}.$$
Apply this inequality to each $z_n:=\frac{1-e(\alpha n-\frac{np}{N})}{e(x-\frac{np}{N})-1}$. We have seen that $|z_n|\lesssim \frac{1}{N\|x-\frac{np}{N}\|}$, and hence
$$\sum_{n\in A}|z_n|^2\lesssim 1.$$
It follows that
$$\prod_{n\in A}\left|1+\frac{1-e(\alpha n-\frac{np}{N})}{e(x-\frac{np}{N})-1}\right|\sim \left|e^{\sum_{n\in A}\frac{1-e(\alpha n-\frac{np}{N})}{e(x-\frac{np}{N})-1}}\right|.$$
Let $\alpha-\frac{p}{N}:=\frac{t}{N^2}$, so $M|t|=1$. Note that since $\|x\|\ge \delta$, it follows that
\begin{equation}
\label{ecc10}
|x-\{\frac{np}{N}\}|<\frac12
\end{equation}
for each  $n\in A$. By invoking Taylor expansions,  \eqref{ecc10}, and  using that
$$|\frac{1}{e(y)-1}-\frac1{2\pi iy}|\lesssim 1$$ for $|y|<\frac12$, we get that
$$\sum_{n\in A}\frac{1-e(\alpha n-\frac{np}{N})}{e(x-\frac{np}{N})-1}=-\sum_{n\in A}\frac{tn}{N^2(x-\{\frac{np}{N}\})}+O(1).$$
We rewrite
$$\sum_{n\in A}\frac{tn}{N^2(x-\{\frac{np}{N}\})}=t\sum_{n=1\atop{\delta\ge|x-\frac{n}{N}|\ge \frac{100}{N}}}^{N}\frac{\frac{n^{*}}{N}}{(Nx-n)},$$
where $n^*:=p^{-1}n$ mod $N$, and $p^{-1}$ is the inverse of $p$ mod $N$.
Our next goal is to prove that
\begin{equation}
\label{ecc12}
\frac{1}{M}\left|\sum_{n=1\atop{\delta\ge|x-\frac{n}{N}|\ge \frac{100}{N}}}^{N}\frac{\frac{n^{*}}{N}}{(Nx-n)}\right|=O(1).
\end{equation}
Since $k$ is odd, it follows from \eqref{ecc4} that $p^{-1}=N_{k-1}$. Let
$$\alpha=\langle a_0,a_1,\ldots\rangle:=a_0+\frac{1}{a_1+\ldots}$$
be the continued fraction expansion of $\alpha$.
We have for each $i\ge 2$
$$p_{i}=a_ip_{i-1}+p_{i-2},$$
$$N_{i}=a_iN_{i-1}+N_{i-2},\,N_0=1,\;N_1=a_1.$$
Due to \eqref{ecc1} we have $a_i\le DM$ for each $i\le k+1$.

Note that $\rho_i:=N_i/N_{i-1}$ satisfies
$$\rho_i=a_i+\frac1\rho_{i-1},\;\rho_1=a_1.$$
Thus,
$$ N/p^{-1}=N_{k}/N_{k-1}=\langle a_k,a_{k-1},\ldots,a_1\rangle.$$
The thing that matters is that all $a_i$ are $O(M)$. Thus, from the recurrence above, the convergents of $N/p^{-1}$, denote them by $M_l/c_l$, have the property that
\begin{equation}
\label{ecc5}
M_{l+1}\lesssim M M_l
\end{equation}
for each $l\le k$ (and similarly for $c_l$, but this will be irrelevant).

It is known that the $l^{th}$ convergent of $p^{-1}/N$ will equal $\frac{c_{l-1}}{M_{l-1}}$, and that the last convergent will equal $p^{-1}/N$. Choose $l_0$ such that
$\frac{N\delta}{M^{3/2}}\lesssim M_{l_0}^{3/2}<{N}\delta$. This is possible due to \eqref{ecc5}. Reasoning as before, we get
$$\frac{1}{M}|\sum_{n=1\atop{|x-\frac{n}{N}|\ge \frac{M_{l_0}^{3/2}}{N}}}^{N}\frac{\frac{n^{*}}{N}}{(Nx-n)}|\lesssim \frac{1}{M}\sum_{N\ge i\gtrsim M_{l_0}^{3/2}}\frac1i\lesssim \frac{\log M+\log(\delta^{-1})}{M}\lesssim_\delta 1.$$

Next, we observe that the remaining part of the sum can be written as
$$\frac{1}{M}\sum_{|j|<M_{l_0}^{3/2}}\frac{\{u+\frac{N_{k-1}j}{N}\}}{j}+O(1),$$
where $u$ is a number whose value is completely irrelevant.

Note that if, say,  $M^5>N$ then the sum above is trivially bounded by $\frac{1}{M}\sum_{|j|<M^5}\frac1{|j|}=O(1)$, and we are fine. Otherwise, we can choose $l_1<l_0$ such that $M^4\lesssim M_{l_1}<M^5$.  The sum above restricted to $|j|\le M_{l_1}^{3/2}$ is trivially $O(1)$.

For $l_1\le l\le l_0-1$ and $M_l^{3/2}\le |j|\le M_{l+1}^{3/2}$, we use that
$$|\frac{N_{k-1}}{N}-\frac{c_l}{M_l}|\le \frac {1}{M_lM_{l+1}},$$
and thus by \eqref{ecc5}
$$|\frac{N_{k-1}j}{N}-\frac{c_lj}{M_l}|\le \frac {M_{l+1}^{3/2}}{M_lM_{l+1}}\lesssim M^{1/2}M_l^{-1/2}.$$
Define
$$C_l:=\{M_l^{3/2}\le |j|\le M_{l+1}^{3/2}: \|u+\frac{c_lj}{M_l}\|\gtrsim M^{1/2}M_l^{-1/2}\}.$$
It follows that
$$|\{M_l^{3/2}\le |j|\le M_{l+1}^{3/2}\}\setminus C_l|\le |\{|j|\le M_{l+1}^{3/2}:\|u+\frac{c_lj}{M_l}\|\lesssim M^{1/2}M_l^{-1/2}\}|$$$$\lesssim M_{l+1}^{3/2}M^{1/2}M_l^{-1/2},$$ and that for each $j\in C_l$
$$|\{u+\frac{N_{k-1}j}{N}\}-\{u+\frac{c_lj}{M_l}\}|\lesssim M^{1/2}M_l^{-1/2}.$$

So we have the following estimate for the error term corresponding to  some $l$
$$\left|\sum_{M_l^{3/2}<|j|<M_{l+1}^{3/2}}\frac{\{u+\frac{N_{k-1}j}{N}\}}{j}-\sum_{M_l^{3/2}<|j|<M_{l+1}^{3/2}}\frac{\{u+\frac{c_lj}{M_l}\}}{j}\right|$$$$\lesssim M^{1/2}M_l^{-1/2}\sum_{j\in C_l}\frac{1}{|j|}+\sum_{M_l^{3/2}\le |j|\le M_{l+1}^{3/2}\atop{j\notin C_l}}\frac{1}{|j|}\lesssim M^2M_l^{-1/2}.$$
Since for each $i$
\begin{equation}
\label{ecc11}
M_{i}\ge M_{i-1}+M_{i-2}\ge 2M_{i-2},
\end{equation}
and since $M_{l_1}\gtrsim M^4$ it follows that the sum of all error terms is bounded by
$$\sum_{l_1\le l}\frac{M^2}{M_l^{1/2}}\lesssim 1$$
as desired. But
$$\sum_{M_l^{3/2}<|j|<M_{l+1}^{3/2}}\frac{\{u+\frac{c_lj}{M_l}\}}{j}=\sum_{r=1}^{M_l}\{u+\frac{c_lr}{M_l}\}\sum_{M_l^{3/2}<|j|<M_{l+1}^{3/2}\atop{j=r \text{ mod }M_l}}\frac{1}{j},$$
and this is $O(\frac{1}{M_l^{1/2}})$, since actually
$$\sup_{P>M_l^{3/2}}|\sum_{M_l^{3/2}<|j|\le P\atop{j=r \text{ mod }M_l}}\frac{1}{j}|=O(\frac{1}{M_l^{3/2}})$$
for each $r$. Summing over $l\ge l_1$ we get using \eqref{ecc11}
$$\sum_{l_0-1\ge l\ge l_1}\left|\sum_{M_l^{3/2}<|k|<M_{l+1}^{3/2}}\frac{\{u+\frac{c_lk}{M_l}\}}{k}\right|\lesssim 1.$$
By putting everything together we conclude that \eqref{ecc12} holds.
\end{proof}

An immediate consequence which only requires trivial modifications is the following.

\begin{corollary}
\label{Coro:745789760598}
Let $A,B\in\C$ with $|A|=|B|=1$. Let also $\alpha$ and $N$ be as in Proposition \ref{P1}. Define $$P(x)=A+Be(\alpha x).$$
Then for each $0<\epsilon<1$ there exist $c_1(\epsilon,A,B,\alpha),c_2(\epsilon,A,B,\alpha)>0$ and a set $P(A,B,\epsilon,\alpha,N)\subset [0,1]$ with measure at least $1-\epsilon$ such that for each $y\in P(A,B,\epsilon,\alpha,N)$
$$c_2(\epsilon,A,B,\alpha)\ge \prod_{n=-N}^{-1}|P(y+n)|\ge c_1(\epsilon,A,B,\alpha)$$
$$c_2(\epsilon,A,B,\alpha)\ge\prod_{n=0}^{N-1}|P(y+n)|\ge c_1(\epsilon,A,B,\alpha).$$

\end{corollary}
The relevance of this result for later applications is that while the sets $P$ are allowed to depend on $N$, the constants $c_1,c_2$ do not depend on $N$.

We can now begin the proof of Theorem \ref{TT2}. By applying the area preserving affine transformations -also called {\em metaplectic transforms}- of the plane (such as translations, rotations, shears, and area one rescalings), it suffices to rule out minimal decay \eqref{eeeeeeeeeeeeeeeeeeeeeee3} for special configurations. See Section 2 in \cite{HRT} for a discussion on this.

Assume for contradiction that there exists a measurable function  $f:\R\to \C$, some $d\in (0,\infty)$ and some $S\subset [0,1]$ with positive measure such that
\begin{equation}
\label{eeeeeeeeeeeeeeeeeeeeeee2}
d<|f(x)|<\infty\;\;\text{for each }x\in S,
\end{equation}
\begin{equation}
\label{eeeeeeeeeeeeeeeeeeeeeee3}
\lim_{|n|\to\infty\atop_{n\in\Z}}f(x+n)=0,
\end{equation}
and
$$f(x+1)(A+Be(\alpha x))=f(x)(E+Fe(\beta x)),$$
for a.e. $x$, for some fixed $A,B,E,F\in\C$, $\alpha,\beta\in\R$, none of them zero. We can also assume $\alpha$ and $\beta$ to be irrational, since the rational case was treated in \cite{De2}. The same metaplectic transforms allow us to assume $0<\alpha<1$.  By re-normalizing, we can trivially assume $E=1$. Let
$$P(x)=A+Be(\alpha x),\;\;\;Q(x)=1+Fe(\beta x).$$
Also, the argument from \cite{De2} shows that the worst case scenario (and the only one that needs to be considered here) is when $|B|=|A|$. Equivalently, $P$ will have zeros. We comment on this in the end of the argument.

By making $S$ a bit smaller, we can also assume that $S+\Z$ contains no zeros of $P$ and $Q$.

Note that by Egoroff's Theorem, \eqref{eeeeeeeeeeeeeeeeeeeeeee3} will allow us to assume (by making $S$ a bit smaller if necessary) that
\begin{equation}
\label{eeeeyurie848eeeee356frsd332}
\lim_{|n|\to\infty\atop_{n\in\Z}}f(x+n)=0,
\end{equation}
uniformly on $S$.

The parameters $D,\alpha,\beta,A,B,F,\epsilon_1,\epsilon_2,\epsilon_3, c_1,c_2,d,m,\gamma$ (some of which are introduced below) will be referred to as {\em fundamental parameters}. They will stay fixed throughout the argument, and in particular will not vary with $N$.

Let us first see how to deal with the contribution coming from the polynomials $Q$. This is done via the {\em conjugates trick} introduced in \cite{De2}. More precisely, let $F=e(\theta)$. Since $S$ has positive measure, it follows that $1_S*1_S$ is continuous and that there exists  an interval $I\subset [0,2]$ and $\epsilon_1>0$ such that
\begin{equation}
\label{ecc14}
1_S*1_S(w)>\epsilon_1
\end{equation}
for each $w\in I$. We can assume without any loss of generality that $I\subset [0,1]$.  There exists $n'\in \N$ large enough such that $m:=[-\frac{2\theta}{\beta}+n'\beta^{-1}]>0$ and $\gamma:=\{-\frac{2\theta}{\beta}+n'\beta^{-1}\}\in I$.  It follows from \eqref{ecc14} that the set $S':=\{x\in S: \gamma-x\in S\}$ has measure at least $\epsilon_1$. The point of this selection is that for each $n\in\Z$, and each $y:=-x-\frac{2\theta}{\beta}+n'\beta^{-1}$, the numbers
$1+Fe(\beta y-n\beta )$ and $1+Fe(\beta x+n\beta)$ are complex conjugates and thus, for each $L\ge 1$, and each $x\in\R$
\begin{equation}
\label{ecc26}
\prod_{n=-L}^{-1}|Q(\gamma-x+n)|=\prod_{n=m+1}^{L+m+1}|Q(x+n)|.
\end{equation}
Let $S''$ be a subset of $S'$ of measure at least $\epsilon_1/2$, and let  $\epsilon_2>0$ depending only on the fundamental parameters $\beta,F$ and $m$ such that
\begin{equation}
\label{ecc16}
\prod_{n=0}^{m}|Q(x+n)|\ge \epsilon_2
\end{equation}
for each $x\in S''$. Let $N$ be as in Corollary \ref{Coro:745789760598}. Let $\epsilon_3>0$ be small enough (depending only on $\epsilon_1$, in particular not depending on $N$) such that the set
$$S(N):=S''\cap\{x\in P(A,B,\epsilon_3,\alpha,N)\}\cap\{x:\gamma-x\in P(A,B,\epsilon_3,\alpha,N)\},$$
has positive measure, and thus is non-empty.  For each $N$ as above, choose a point $x_N\in S(N)$. Let $z_N:=\gamma-x_N$. The recurrence along the orbits of $x_N$ and $z_N$ implies that
$$|f(x_N+N+m+2)|=|f(x_N)|\frac{\prod_{n=0}^{N+m+1}|Q(x_N+n)|}{\prod_{n=0}^{N+m+1}|P(x_N+n)|}$$
$$|f(z_N-N)|=|f(z_N)|\frac{\prod_{-N}^{n=-1}|P(z_N+n)|}{\prod_{-N}^{n=-1}|Q(z_N+n)|}.$$
Multiply these equalities. Using the fact that $x_N,z_N$ are in $S$, \eqref{eeeeeeeeeeeeeeeeeeeeeee2}, \eqref{ecc26} with $x:=x_N$ and $L:=N$, \eqref{ecc16} with $x:=x_N$, Corollary \ref{Coro:745789760598} and the fact that
$$\prod_{n=N}^{N+m+1}|P(x_N+n)|\le (2|A|)^{m+2},$$
it follows that
$$|f(x_N+N+m+2)||f(z_N-N)|\ge \frac{d^2\epsilon_2c_1(\epsilon_3,A,B,\alpha)}{(2|A|)^{m+2}c_2(\epsilon_3,A,B,\alpha)}.$$
The important thing is that the constant on the right depends only on the fundamental parameters, and not on $N$. By letting $N\to\infty$, this will contradict the uniformity assumption \eqref{eeeeyurie848eeeee356frsd332}. This ends the proof of Theorem \ref{TT2}, under the assumption that $|A|=|B|$.

If $|A|\not=|B|$, then things are much easier, and have already been addressed in \cite{De2}. We briefly recap the argument. By invoking Riemann sums and the fact that the derivative of $\phi(x):=\ln|A+Be(x)|$ satisfies
$$\inf_{x\in [0,1]}|\phi'(x)|\gtrsim_{A,B} 1,$$
we get that
$$|\sum_{n=0}^{N-1}\ln |P(x+n)|-N\int_{0}^1\phi|\lesssim_{A,B} 1$$
$$|\sum_{n=-N}^{-1}\ln |P(x+n)|-N\int_{0}^1\phi|\lesssim_{A,B} 1,$$
for each $x\in [0,1]$ and each $N$ such that
$$N\|N\alpha\|\le 1.$$
In particular,

$$|\sum_{n=0}^{N-1}\ln |P(x_N+n)|-\sum_{n=-N}^{-1}\ln |P(z_N+n)||\lesssim_{A,B} 1$$
and thus
$$\frac{\prod_{-N}^{n=-1}|P(z_N+n)|}{\prod_{n=0}^{N-1}|P(x_N+n)|}\sim_{A,B} 1,.$$
This will replace Corollary \ref{Coro:745789760598} in the argument above. Everything else will be the same.

\end{document}